\documentclass [winedt,yap]{iitparc}
\usepackage{cite}
\usepackage{amsmath,amssymb,amsfonts,bm}
\usepackage{lscape}
\usepackage{floatfig,wrapfig,epsfig}
\usepackage{subfigure}
\usepackage{color}
\usepackage{rotating}
\usepackage{curves}
\usepackage{ulem}
\usepackage[mathscr]{eucal}
\usepackage{array}            
\RequirePackage{psfrag}
\RequirePackage{graphicx}

                \usepackage{accents}
                \DeclareMathSymbol{\widehatsym}{\mathord}{largesymbols}{"62}
                \newcommand\lowerwidehatsym{%
                  \text{\smash{\raisebox{-1.3ex}{%
                    $\widehatsym$}}}}
                \newcommand\fhat[1]{%
                  \mathchoice
                    {\accentset{\displaystyle\lowerwidehatsym}{#1}}
                    {\accentset{\textstyle\lowerwidehatsym}{#1}}
                    {\accentset{\scriptstyle\lowerwidehatsym}{#1}}
                    {\accentset{\scriptscriptstyle\lowerwidehatsym}{#1}}
                }

\begin{document}\normalem
\initfloatingfigs
\frontmatter          

\IssuePrice{25.00}%
\TransYearOfIssue{2018}%
\TransCopyrightYear{2018}%
\OrigYearOfIssue{2016}%
\OrigCopyrightYear{2016}%

\TransVolumeNo{79}%
\TransIssueNo{4}
\OrigIssueNo{62}%
\OrigPages{169-–187} 
\OrigCopyrightedAuthors{%
{P.Yu.~Chebotarev},
{V.A.~Malyshev},
{Ya.Yu.~Tsodikova},
{A.K.~Loginov},
{Z.M.~Lezina},
{V.A.~Afonkin}%
}
\OrigJournalName{Upravlenie Bol’shimi Sistemami}

\mainmatter

\setcounter{page}{500}

\Rubrika{LARGE SCALE SYSTEMS CONTROL}


\def\x#1{} 
\def\xz{\hspace{-.17em}}
\def\MP{M\xz P}
\def\ell{n}
\def\a{\alpha}
\def\Up#1{\vspace{-#1em}}
\def\cdc{,\ldots,}
\def\xz{\hspace{-.07em}}
\def\xy{\hspace{.07em}}
\def\beq#1{\begin{equation}\label{#1}}
\def\eeq{\end{equation}}
\def\M{\fhat{\rm M}}                             

\title{The Optimal Majority Threshold as a~Function of the~Variation Coefficient of the~Environment}
\author{%
{P.~Yu.~Chebotarev}$^{\xy*,**,***a}$, 
{V.~A.~Malyshev}$^{\xy*,**,***,b}$,
{Ya.~Yu.~Tsodikova}$^{\xy*,c}$,
{A.~K.~Loginov}$^{\xy*,d}$,
{Z.~M.~Lezina}$^{\xy*,e}$,
{V.~A.~Afonkin}$^{\xy*,**,f}$\\
}
\institute{$^*$Trapeznikov Institute of Control Sciences, Russian Academy of Sciences, Moscow, Russia\\
           $^{**}$Moscow Institute of Physics and Technology, Moscow, Russia\\
           $^{***}$Kotelnikov Institute of Radioengineering and Electronics, Russian Academy of Sciences, Moscow, Russia\\
           e-mail: $^a$pavel4e@gmail.com, $^b$vit312@gmail.com, $^c$codikova@mail.ru, $^d$a\_k\_log@mail.ru, $^e$lezinazo@gmail.com, $^f$afonkinvadim@yandex.ru}

\received{Received January 19, 2016}
\authorrunning{CHEBOTAREV, MALYSHEV, TSODIKOVA, LOGINOV, LEZINA, AFONKIN}
\titlerunning{OPTIMAL MAJORITY THRESHOLD}

\maketitle

\begin{abstract}
Within the model of social dynamics determined by collective decisions in a stochastic environment $($ViSE model$)$, we consider the case of a homogeneous society consisting of classically rational economic agents (or homines economici, or egoists). We present expressions for the optimal majority threshold and the maximum expected capital increment as functions of the parameters of the environment. An estimate of the rate of change of the optimal threshold at zero is given, which is an absolute constant: $(\sqrt{2/\pi}-\sqrt{\pi/2})/2$.


\medskip\noindent$\!\!\!\!\!\!\!\!$\DOI{xxxx} 
\end{abstract}

\setcounter{footnote}{1}
\begin{flushright}\x{2}
{\emph{An egoist is like someone sitting in a well since long ago.}\x{"}

Kozma Prutkov, ``Fruits of Reflection'' (1853--1854)}
\end{flushright}

\section{Introduction}
\label{s_intro}

Consider the ViSE (Voting in Stochastic Environment) model \cite{Che2706AiT} in the case where the \emph{society\/} consists of $n$ {classically rational economic agents\/} (\emph{homines economici} \cite{OBoyle09}), who are {boundedly rational egoists\/} (hereafter, \emph{egoists}). Such an agent maximizes his/her individual utility function in every act of choice. This is obviously the most profitable \emph{individual\/} strategy. Cooperative and altruistic strategies within the ViSE model have been considered in \cite{Che2706AiT,Che06AiT,Che08PU,CLCLB09IPU}.

Let $\a\in[0,1]$ be a strict relative \emph{voting threshold\/}, which means that any proposal is accepted and implemented if and only if the proportion of the society supporting this proposal is greater than~$\alpha.$

Each \emph{participant} (\emph{agent\/}) is characterized by the current value of \emph{capital} (which can also be interpreted as the value of individual \emph{utility}). \emph{A proposal} [\emph{of the environment\/}] is a vector of proposed capital increments of the participants. This concept allows to model potential innovations that are beneficial for some agents and disadvantageous for others. As a result of the implementation of such a proposal, the capitals of the agents of the first type increase, while the capitals of the remaining agents decrease.

The proposals are consequently put to a general vote. Each homo economicus votes for those and only those proposals that increase his/her individual capital (utility). If a proposal is supported by a proportion of the society exceeding the threshold $\a$, then it is accepted (the voting procedure is ``$\a$-majority'' \cite{NitzanParoush1982optimal,NitzanParoush84qualified,Rae69APSR-Decision}) and the participants' capitals get the proposed increments. Otherwise, all capitals remain the same. 
The voting threshold $\a$ will also be called \emph{majority threshold\/} and, more precisely, \emph{acceptance threshold\/}, since $\a<0.5$ is allowed.

The generation and adoption of proposals are repeated over and over again, whereas the subject of study is the change of the participants' capital as a result of this process. Does it inevitably lead to an increase in the social welfare or can democratic decisions systematically reduce the total capital of the society? Does the financial inequality grow? How many participants are ruined?

In accordance with the basic ViSE model, the capital increments that form the proposals of the environment are the realizations of independent identically distributed random variables. In this paper, we study the case where these variables have the normal distribution $N(\mu,\sigma)$ whose mathematical expectation is $\mu$ and standard deviation is~$\sigma.$%

The ratio~$\sigma/\mu$ is called the \emph{coefficient of variation\/} of a random variable. In what follows, we need the \emph{inverse coefficient of variation\/}: $\rho=\mu/\sigma,$ which we call the \emph{adjusted\/} or \emph{normalized mean of the environment\/}.
If $\rho>0,$ then the opportunities provided by the environment are favorable on average; if $\rho<0,$ then the environment is unfavorable.

An aspect of social practice that can be examined using the ViSE model is adoption by the Parliament of various bills that are ``prompted by life'' (environment), i.e., by economic and/or political conjuncture. It can be assumed that the members of Parliament, being lobbyists of certain interests and adherents of certain beliefs, are so interested in accepting or rejecting a bill (e.g., a budget draft) that this concernment is adequately expressed in terms of individual utility or capital. Of course, the adoption of any other collective decisions can also be considered (to some extent) in terms of the ViSE model and its modifications.

In the present paper, the following topics are studied:

\Up{.8}
\begin{itemize}
\item the dynamics of the participants' capital under the stated model assumptions;
\item the \emph{optimal acceptance threshold\/}, i.e., the threshold that maximizes the total capital of the society;
\item dependence of the optimal acceptance threshold on the parameter $\rho=\mu/\sigma$ characterizing the favorability of the environment.
\end{itemize}

\section{Dependence of capital increments on the number of agents and the parameters of the environment}
\label{s_onEnv}

The dependence of the average (mean, expected) one-step capital increment on the parameters of the environment and the acceptance threshold has been studied in~\cite{CheLogTsoLez13HSE8}.
In particular, it was found that voting by a simple majority (i.e., with $\a=0{.}5$) in a moderately unfavorable environment leads to a decrease of the total capital.

Analytically, this dependence can be expressed by the following proposition.
\begin{proposition}\label{p_Md}
$1.$ In a society of $\ell$ egoists$,$ the \emph{expected\/} one-step capital increment $(\tilde d)$ of an agent is\/$:$
\beq{onlyEgoExp}
     {\rm M}(\tilde d) 
                              \,=\,\sigma\!\!\sum_{x=[\a \ell]+1}^\ell\left(\rho+\tfrac f q\!\left(\tfrac x{p\ell}-1\right)\right)b(x|\ell),
\eeq
where $[y]$ is the integer part of $y,$
$\,p = F(\rho),$
$\,q = 1 - p = F(-\rho),$ $f=f(\rho),$
$\,b(x\,|\,\ell) = \begin{pmatrix}\ell\\x\end{pmatrix} p^x q^{\ell-x},$
$\,F(\cdot)$ and $f(\cdot)$ being the standard normal distribution function and the standard normal density\/$,$ respectively.

$2.$ The standard normal approximation of the binomial distribution leads to the approximate formula
\beq{onlyEgoExp1}
\M(\tilde d) = \sigma\!\left(\rho F(\tau) + \tfrac f{\sqrt{qp \ell}}f(\tau)\right),
\eeq
where
\beq{e_tau}
\tau = \frac{p\ell-[\a\ell]-0{.}5}{\sqrt{qp\ell}}.
\eeq
\end{proposition}

Proposition~\ref{p_Md} follows from the lemmas on the ``normal voting sample''~\cite{Che06AiT}.
It can be observed that $p$ and $q$ are the probabilities of positivity and negativity for a single capital increment, respectively, while $f/\sigma$ is the density of a zero increment in a proposal.

A normal approximation for the binomial distribution is recommended when $qpn\ge9$. For a fixed $qpn$, its accuracy is maximal for $p=0{.}5$ and decreases when $p$ approaches $0$ or~$1$. That is why for $0{.}1<p<0{.}9$, the normal approximation is frequently adopted whenever $qpn>5$. For $p$ very close to $0$ or $1$, sometimes $qpn>25$ is required.

\begin{figure}[ht!]
\begin{center}
\includegraphics[width=3.6in]{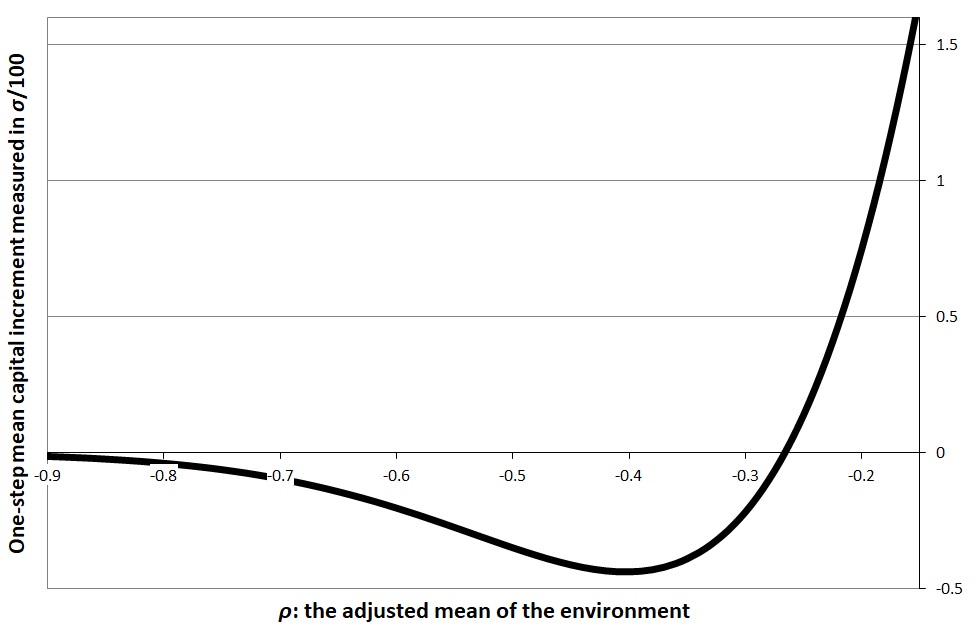} 
\end{center}\Up{1.1}
\caption{\label{f_}One-step mean capital increment of an agent$:$ $21$ agents; $\a= 0{.}5.$}
\end{figure}

The dependence of the one-step mean capital increment ${\rm M}(\tilde d)$ of an agent on the adjusted mean of the environment $\rho$ for $21$ participants and $\a=0{.}5$ is shown in Fig.\:\ref{f_}.

Figure \:\ref{f_} shows that for $\rho\in(-0{.}9,-0{.}266)$, the mean capital increment is an appreciable negative value, i.e., proposals approved by the majority are, on average, unprofitable (and ``confiscatory'') for the society. The corresponding part of the curve will be called the ``\emph{pit of losses}\/.'' For $\rho<-0{.}9$, the negative ${\rm M}(\tilde d)$ is very close to zero, since the proposals are extremely rarely accepted.

The phenomenon of ``pit of losses'' has much in common 
with the results of A.V.\:Malishevsky (described in \cite[Chapter\;2, Section\;1.3]{Mirkin74}), 
who was the first to demonstrate (in 1969) that the eventual outcome of a series of egoists' votes can be extremely unprofitable for all of them. It also evokes association with the aphorism of Kozma Prutkov\footnote{Gibian, G., editor, \emph{The Portable Nineteenth-century Russian Reader}, New York: Penguin Books, 1993.}, ``An egoist is like someone sitting in a well since long ago,'' used as an epigraph to this paper.

This phenomenon can be explained by the fact that, due to the negativity of $\rho,$ positive increments proposed by the environment have, on average, smaller absolute values than negative ones. As a result, the total loss of the losing minority systematically exceeds the total income of the winning majority. Hence, despite the desire of all agents to increase their capital and majority approval of all the decisions, social welfare reduces.
Thus, in a moderately unfavorable environment, decisions taken by a simple majority of egoists reduce their total capital!

\begin{figure}[ht!]
\begin{center}
\includegraphics[width=3.8in]{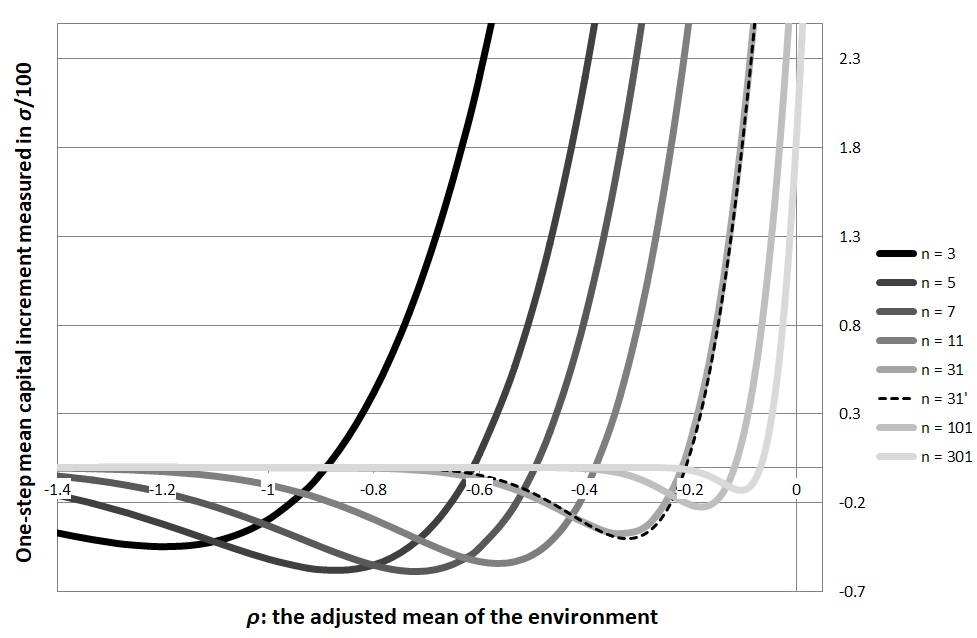} 
\end{center}\Up{1.1}
\caption{\label{f_n}One-step mean capital increment of an agent for odd $n$'s and $\a= 0{.}5.$}
\end{figure}

How does the mean capital increment depend on the number of agents? This dependence is shown in Fig.~\ref{f_n}.
With an increase of $n$, the minimum point moves to the right and the ``pit of losses'' becomes narrower. The depth of the ``pit'' increases at first, reaches a maximum at $n=7,$ and then tends to zero monotonically with further increase of~$n$. Let us note that for small $n$'s, the approximation of \eqref{onlyEgoExp1} has insufficient accuracy and does not reflect the initial tendency of the ``pit'' to deepen with the growth of~$n.$ For $n=31,$ its accuracy (see the dashed line in Fig.\:\ref{f_n}) becomes acceptable.

Thus, the effect of capital reduction caused by majority decisions gradually weakens with the growth of~$n$. For even $n$'s the ``pit'' is shallower than for the neighboring odd numbers of agents, since the \emph{de facto\/} majority threshold is higher in the latter case: to accept a proposal, an overbalance in at least two votes is necessary.

It can be observed that the shape of the ${\rm M}(\tilde d)$ curves is similar for different odd~$n$'s. Moreover, if we denote the dependence of ${\rm M}(\tilde d)$ on $\rho$ for $n$ agents and a fixed $\sigma$ by $\varphi_n(\rho),$ then with $n$ and $n'$ exceeding $20$ we have a quite accurate approximate relationship:
\beq{e_appn}
\varphi_{n'}(\rho) \approx \sqrt{\tfrac n{n'}}\,\varphi_n\Big(\rho\,\sqrt{\tfrac{n'}n}\,\Big).
\eeq

\begin{figure}[h!]
\begin{center}
\includegraphics[width=3.6in]{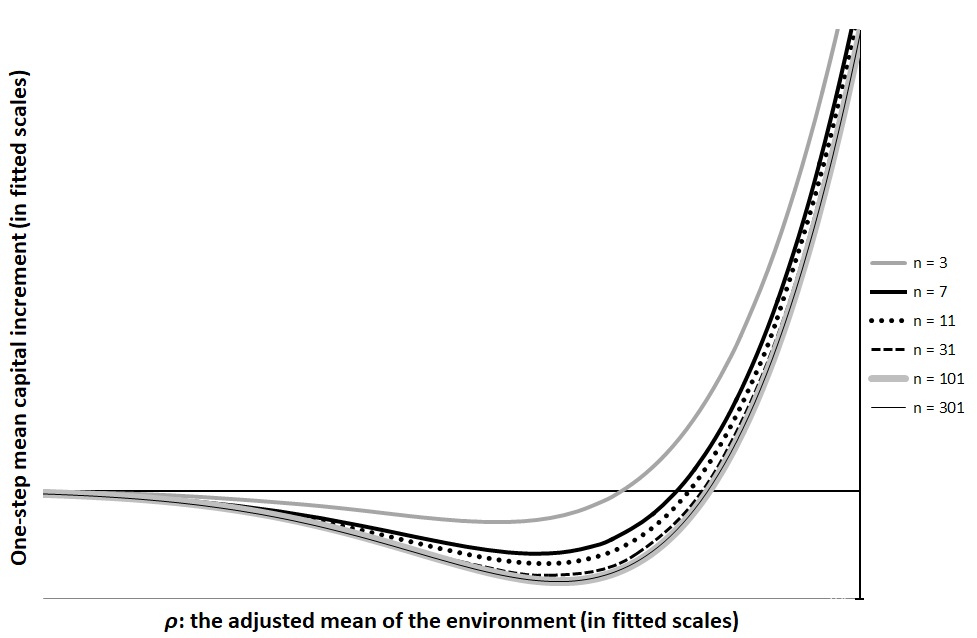} 
\end{center}\Up{1.1}
\caption{\label{f_c}Matching curves of Fig.\:$\ref{f_n}$ by means of correspondence~\eqref{e_appn}.} 
\end{figure}

Thus, a $k^2$-fold decrease in $n$ leads to a $k$-fold stretching of the curve along both axes.
A matching of the curves based on \eqref{e_appn} is illustrated by Fig.~\ref{f_c}, where the functions $\sqrt n\,\varphi_n\big(\rho\,/\sqrt n\,\big)$ with various $n$ are depicted.
If $n$ and $n'$ are greater than 30, then the accuracy of \eqref{e_appn} is quite high.

It follows from \eqref{onlyEgoExp1} and \eqref{e_appn} that a $k^2$-fold increase in the number of participants has a similar effect as a
$k$-fold decrease in~$\sigma$. Hence, an increase in variance ``balances'' the same increase in the number of voters: the graph of the dependence of ${\rm M}(\tilde d)$ on $\rho$ remains almost unchanged under these simultaneous transformations.

Finally, observe that for an odd $n,$ $\a=0{.}5,$ and $\rho=0$ (neutral environment), \eqref{onlyEgoExp1} and \eqref{e_tau} yield a very simple formula:
\beq{e_beauty1}
\M(\tilde d)
=\frac\sigma{\pi\sqrt n}.
\eeq

On the other hand, such a result of averaging with rejection within a probabilistic model looks quite natural. %

\section{Dependence of the capital dynamics on the majority threshold} 
\label{s_onThr}

\begin{figure}[ht!]
\begin{center}
\includegraphics[width=3.7in]{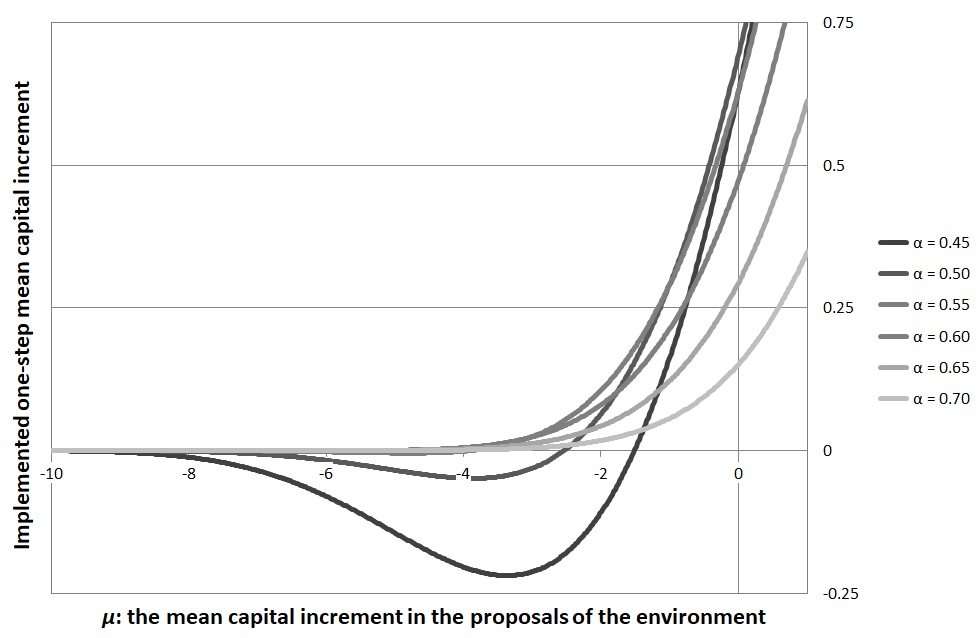} 
\end{center}\Up{1.1}
\caption{\label{f_a}One-step mean capital increment of an agent for various acceptance thresholds $\a\; (n=21,\; \sigma=10).$}
\end{figure}

As one could guess, the ``pit of losses'' becomes narrower and shallower with an increase in the majority threshold~$\a$: the mean capital increment of an agent falls less in the negative area (Fig.\:\ref{f_a}).

At the same time, looking to the right we notice that the higher $\a$ is, the slower the agent's capital grows with an increase of~$\mu.$ 
The society can insure itself against loss under a moderately negative $\mu$ by choosing a high enough $\a$, however, this reduces the benefit of a large $\mu$ when the environment becomes favorable.

Fig.\:\ref{f_a} (where approximation \eqref{onlyEgoExp1} is used) contains, among others, the curve for $\a=0{.}45$. It has a relatively wide and deep ``pit of losses,'' however, it grows faster than the curves with $\a>0{.}45$ as $\mu$ increases. Thereby it surpasses the curves corresponding to $\a\ge0{.}5$ in the domain of positive $\mu$ (cf.\ Fig.\:\ref{f_a20},{\it a}). This means that in a favorable environment, it is rational to ``take risks" of accepting proposals that are only supported by a sufficiently large minority. 
Despite the fact that these can be disadvantageous for the majority, the benefit of the supporting minority will systematically exceed the loss of the majority. After a sufficient number of steps, this will lead to profits for all, since the agents' capital increments are mutually independent.

\section{Optimum acceptance threshold in a society of egoists}
\label{s_optthr}

The above analysis shows that for any $n,\,$ $\sigma,$ and the \emph{environment favorability\/} $\mu,$ there is a kind of \emph{optimum\/\footnote{On other approaches to optimizing the majority threshold see
\cite{AzrieliKim14Pareto, NitzanParoush1982optimal} and \cite{Rae69APSR-Decision,SekiguchiOhtsuki15} for the case of multiple voting.} acceptance threshold\/}~$\a$. It is $\a$ that maximizes the average capital increment of the society.
\begin{figure}[ht!]
\begin{center}
{\small ({a})\hspace{24.2em} ({b})}\\
\includegraphics[width=3.44in]{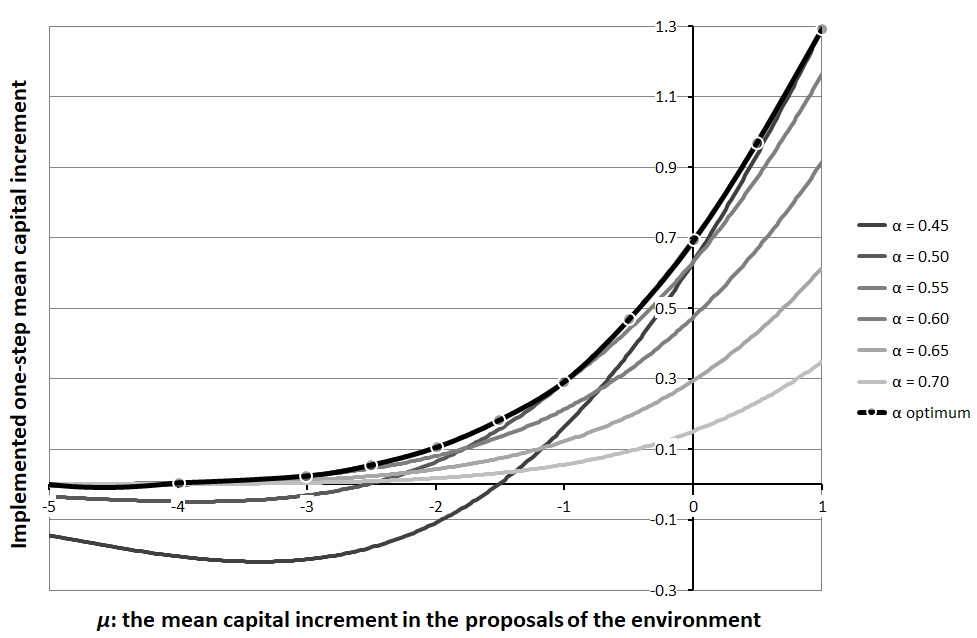}
\includegraphics[width=3.44in]{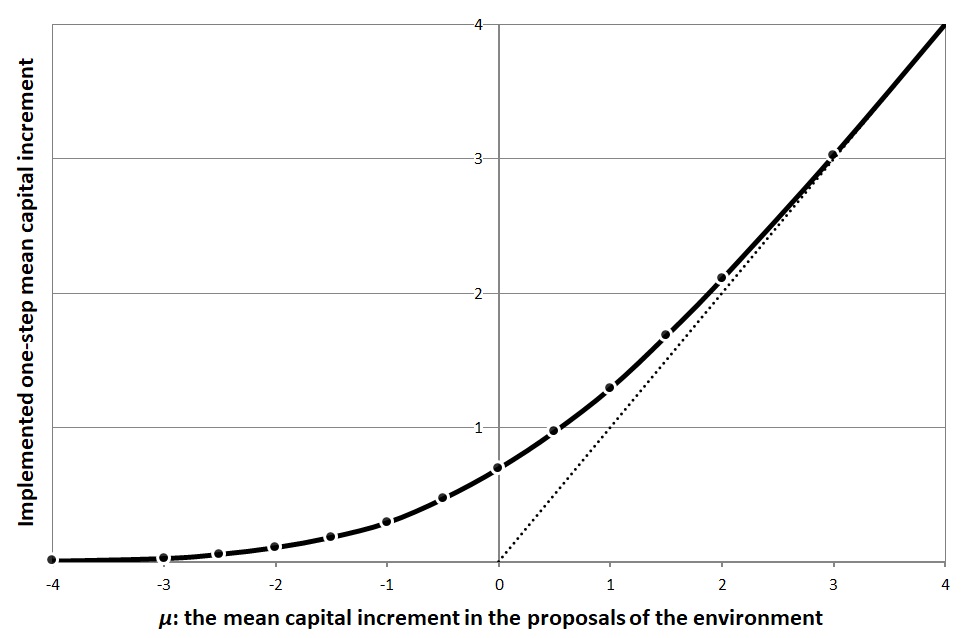}
\end{center}\Up{1.1}
\caption{\label{f_a20} The spline that maximizes the increment of capital w.r.t.~$\a$ {(a)} and
                       the average capital increment of an agent in voting with the optimal acceptance thresholds {(b)}; $n=21,\; \sigma=10.$}
\end{figure}

The optimal $\a$ can be found graphically. For each $\mu_0,$ among all different curves of ${\rm M}(\tilde d)$ versus $\mu$ corresponding to various $\a$ (like those in Fig.\;\ref{f_a}) find the one having the largest ordinate at abscissa~$\mu_0.$ This curve matches optimal~$\a$'s for~$\mu_0.$
Now we can construct the curve of ${\rm M}(\tilde d)$ versus $\mu$ for which an optimal $\a$ is taken at each $\mu.$
It is the spline of the ``leading" fragments of the curves such as those shown in Fig.\:\ref{f_a} (see Fig.\:\ref{f_a20}). Note that it is not the upper envelope\footnote{An envelope touches all curves in an \emph{infinite\/} set.} of the set of curves, since only $n+2$ essentially different voting thresholds~$\a$ are possible for each $n$: $\{-\tfrac1n, 0, \tfrac1n, \tfrac2n\cdc1\}$ (the first and the last values correspond to the acceptance and rejection of all proposals, respectively).
The constructed curve provides the \emph{maximum\/} (w.r.t.~$\a$) \emph{possible mean increment\/} of agent's capital for each~$\mu.$

It is easy to see (Fig.\:\ref{f_a20},{\it b}) that voting with optimal acceptance thresholds is devoid of ``pits of losses'' and always yields positive expected capital increments.

\begin{figure}[h!]
\begin{center}
\includegraphics[width=3.9in]{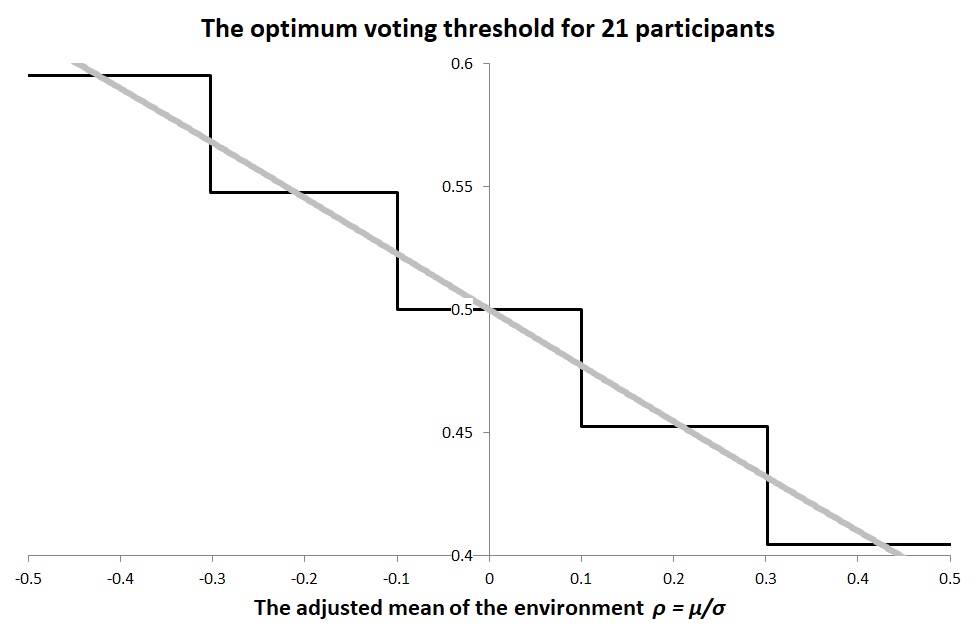} 
\end{center}\Up{1.5}
\caption{\label{f_as}Mean values of the equivalence classes of optimal acceptance thresholds~$\a_0$ $($a ``ladder''$)$ for $n=21$
along with the approximation \eqref{optAlpha} of the optimal acceptance threshold $($the gray line$)$.}
\end{figure}

An expression for the optimal acceptance threshold will be obtained in Section~\ref{s_opthrF}. Now observe that its dependence on $\mu$ is a ``ladder'' with steps of equal height due to the aforementioned finiteness of the set of essentially different thresholds.
If $\a$ is an optimal threshold and $[\a_1 n]=[\a n],$ then $\a_1$ is also optimal.
Fig.\:\ref{f_as} presents the graph of the mean values of the equivalence classes of optimal acceptance thresholds as a function of $\rho\in[-0{.}5, 0{.}5]$ for $n=21$; vertical lines are added for clarity.

If $\bar\a_0$ is the mean value of the equivalence class of optimal thresholds for fixed $n,\sigma$, and $\mu$, then this class is the half-interval $[\bar\a_0-\tfrac1{2n}, \bar\a_0+\tfrac1{2n}[\xy.$
Outside the segment $\rho\in[-0{.}7,\, 0{.}7]$, if the acceptance threshold $\a$ is close to the optimal one and the number of agents is appreciable, then proposals are almost always accepted (to the right of the segment) or hardly ever accepted (to the left of the segment) (cf.\ Fig.\:\ref{f_a20},{\it b}). Therefore, the issue of determining the exact optimal threshold loses its practical value in this case.

Finding the optimal acceptance threshold $\a_0$ in real-world situations looks like a solvable problem. To estimate $\a_0$, it is sufficient to estimate $\rho=\mu/\sigma$ on the basis of statistics and to have reason to believe that the ViSE model is relevant for the situation under consideration at least approximately.
However, even if the estimation of $\rho$ and the question of ViSE model adequacy cause difficulties, the general conclusion that the increase of the acceptance threshold is reasonable when the environment becomes less favorable, seems to remain true.
This conclusion is based on the fundamental fact that a total loss of a minority in an unfavorable environment can systematically exceed the total income of the majority supporting a proposal. 
This fact determines a common flaw of standard voting procedures: the votes are taken into account irrespective of the importance of the considered issue for each voter; the total loss/profit caused by a proposal is ignored. 

Now we turn to derive an approximate expression for the optimal acceptance threshold in a society consisting of homines economici.

\section{Expressions for the optimal acceptance threshold}
\label{s_opthrF}

The following theorem gives an approximation of the optimal voting threshold~$\a_0$. %
\begin{theorem}
\label{pro:0}
Let $\a_0$ be the acceptance threshold that maximizes the expected one-step capital increment in the society of egoists. 
Then the standard normal approximation of the binomial distribution leads to\xy$:$\\
\Up{1.7}
\begin{itemize}
\item the estimate 
\beq{optAlpha}
\hat\a_0 = p\left(1 - \frac{q\rho}f\right)\xz;
\eeq
\item the expression 
\beq{onlyEgoExp2}
\M(\tilde d_0) = \mu\xy F\xz\left(\frac\mu\nu\right) + \nu f\xz\left(\frac\mu\nu\right)
\eeq
for the corresponding maximum expected capital increment\/$,$ where $\nu=\sigma\xz f/\sqrt{qp \ell}$ and the remaining notations are given in Section~$\ref{s_intro}$ and Proposition~$\ref{p_Md}.$
\end{itemize}
\end{theorem}

The proofs are given in the Appendix.

Note that threshold~\eqref{optAlpha} is solely determined by the environment parameters and does not depend on~$\ell$. Since $p,\,q,$ and $f$ are functions of $\rho,$ the only parameter that determines~$\hat\a_0$ is~$\rho.$ The detailed formula is:
$$
\hat\a_0 = F(\rho)\left(1 - \frac{\rho\xy F(-\rho)}{f(\rho)}\right).
$$

The ``ladder'' dependence shown in Fig.\:\ref{f_as} can be obtained from~\eqref{optAlpha} by applying transformation
\beq{e_ladd}
\bar\a_0 = \frac{[\hat\a_0n]+0{.}5}n.
\eeq

Theorem~\ref{pro:0} provides extremely simple analytical approximations for the optimal threshold and the corresponding expected capital increment of an agent. For comparison, the problem of finding the abscissas and ordinates of the functions ${\rm M}(\tilde d)$ minima (shown in Figures~\ref{f_n} and~\ref{f_a}) and the points of their intersection with the horizontal axis does not lead to simple analytical expressions.

\section{Rate of change of the optimal voting threshold as a function of~$\bm\rho$}
\label{s_Velo}

\begin{figure}[h]
\begin{center}
\includegraphics[width=3.8in]{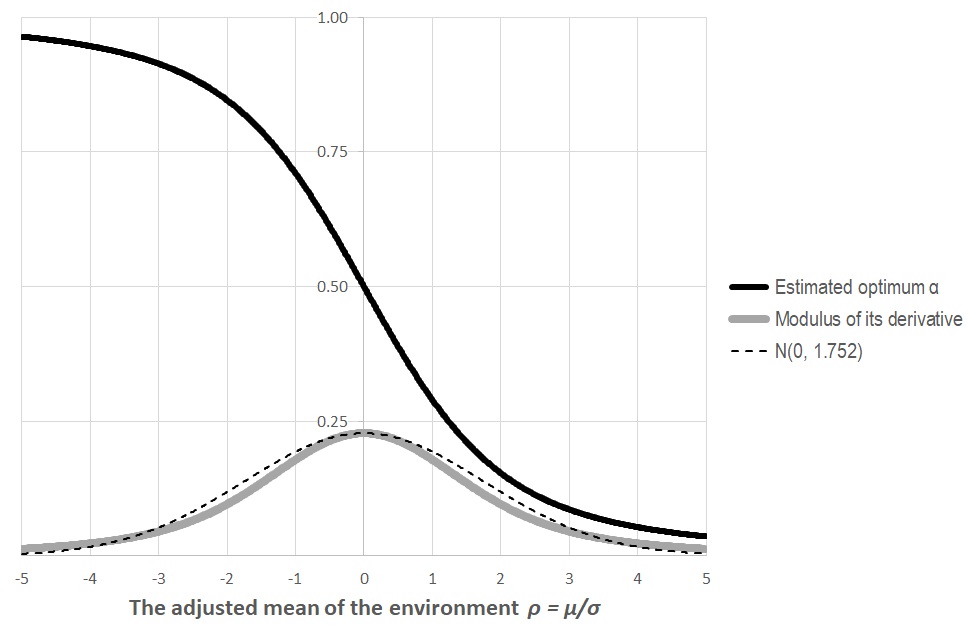} 
\end{center}\Up{1.1}
\caption{\label{f_55}Approximation $\hat\a_0$ of the optimal acceptance threshold and the derivative of the expression~\eqref{optAlpha} taken with minus. The latter interpreted as a distribution density has ``relatively heavy tails,'' which is noticeable when comparing it with the normal density indicated by a dotted line.}
\end{figure}

The dependence of $\hat\a_0$ on $\rho$ in Fig.~\ref{f_as} looks linear, however, linearity is inevitably violated on a wider range of arguments due to the finiteness of $\a$ and the infiniteness of $\rho$ (Fig.~\ref{f_55}).

To find the slope of the curve $\hat\a_0(\rho)$ at zero, we differentiate the function~\eqref{optAlpha} with respect to~$\rho.$
\begin{proposition}\label{p_dera}
\beq{e_dera}
\frac{d\hat\a_0(\rho)}{d\rho} = \frac{(f+p\rho)(f-q\rho)-qp}{f}.
\eeq
\end{proposition}

The function ${d\hat\a_0(\rho)}/{d\rho}$ is negative; its absolute value is shown in Fig.\:\ref{f_55}.

In the case of neutral environment ($\rho=0$), we have $p=q=1/2$ and $f=1/\sqrt{2\pi}.$ Substituting these into~\eqref{e_dera} yields
\begin{corollary}\label{co_limlim} %
$$
\frac{d\hat\a_0(\rho)}{d\rho}\Big|_{\rho=0} = \frac12\left(\sqrt\frac2\pi-\sqrt\frac\pi2\right) \approx -0{.}2277.
$$
\end{corollary}

This derivative at zero does not depend on either the number of agents or the other parameters of the model. We have: if $\mu$ decreases from $0$ to $-\tfrac\sigma2$ $($i.e., $\rho$ changed from $0$ to $-\tfrac12)$ due, say, to an economic crisis, then the optimal majority threshold increases from 50\% to 61\%. If $\mu$ decreases to $-\sigma,$ then 
$\hat\a_0\approx 71\%.$

\section{Conclusion}

In this paper, we have obtained and interpreted several relationships that describe voting with an optimal (i.e., maximizing the total capital of the society) majority threshold in the assumptions of the ViSE model. These relationships led us to the conclusion that the acceptance threshold should be increased if the environment generating proposals is becoming less favorable and should be decreased in the opposite nonstable case. The optimal acceptance threshold has been estimated by means of a function of the environment parameters, which does not depend on the number of participants. The corresponding maximum expected capital increment has been specified.

We have also obtained analytical expressions for the ``pit of losses'', i.e., the graph of the function representing the rate of decrease of egoists' capital/utility caused by their democratic decisions in a moderately unfavorable environment.

The derivative at zero of the estimated optimal acceptance threshold with respect to the environment favorability is a constant $\big(\sqrt{2/\pi}-\sqrt{\pi/2}\big)/2$.
The results of this paper allow interpretation in terms of making real-world collective decisions in the environment whose favorability and variability can be empirically assessed.

\section*{Acknowledgments}
This work was supported by the Russian Science Foundation, project 16-11-00063 granted to IRE~RAS.

\appendix{}

\PTH{\ref{pro:0}}
To find the argument of the maximum value of the capital increment by differentiation, we replace $[\a\ell]+0{.}5$ by $\a n$ in the expression~\eqref{e_tau} (cf.~\eqref{e_ladd}). With this substitution, we obtain a differentiable function that intersects the original one at points $\Big\{\frac{k+0{.}5}n\,|\,k=0\cdc\ell\Big\}$. Now instead of~\eqref{e_tau}, we have
\beq{e_tau1}
\tau = (p-\a)\sqrt{\frac \ell{qp}}\,.
\eeq

Differentiating expression \eqref{onlyEgoExp1} with respect to $\a$ and substituting~\eqref{e_tau1} (the result of substitution is denoted by $\M(\a,\tau,\sigma)$), we find:
\beq{derivM}
    \frac{d\M(\a,\tau,\sigma)}{d \a}
    = \mu f(\tau) \frac{d \tau}{d \a} + \frac{\sigma f}{\sqrt{qp\ell}}f(\tau)\xy\tau\frac{d (-\tau)}{d \a}.
\eeq

The first-order condition for the maximum (the first derivative being zero) is reduced to
\beq{firstOrderCondition}
    -\rho = \frac{(\a-p)f}{qp}.
\eeq
Taking into account the negativity of the second derivative, we obtain~\eqref{optAlpha}.

The second assertion of the theorem is proved by substituting the approximation \eqref{optAlpha} of the optimal acceptance threshold into~\eqref{onlyEgoExp1}.
Theorem~\ref{pro:0} is proved.

\PPR{\ref{p_dera}}
Taking into account that $dp/d\rho=-dq/d\rho=f$ and $df/d\rho=-\rho f$ \Big(since $f=\tfrac1{\sqrt{2\pi}}e^{-\rho^2/2}$\Big), we obtain
\begin{eqnarray*}
\frac{d\hat\a_0(\rho)}{d\rho}
&\!=\!& \frac d{d\rho}\:p\xz\left(1-\frac{q\rho}f\right)\\ 
&\!=\!& f\xz\left(1-\frac{q\rho}f\right) + p\xz\left(-\frac{(-f)\rho f+qf+\rho fq\rho}{f^2}\right)\\ 
&\!=\!& \frac{(f+p\rho)(f-q\rho)-qp}{f}.
\end{eqnarray*}
Proposition~\ref{p_dera} is proved.




\end{document}